\documentclass[a4paper]{amsart} 

\usepackage[all]{xypic} 
\usepackage{amssymb,amsmath}

\def\p{{\mathbb P}}
\def\zz{{\mathbb Z}}
\def\oc{{\mathcal O}}
\def\su{\subset} 
\def\op1{{\mathcal O}_{\p^1}}

\def\q{\mathcal Q}
\def\j{\mathcal J}
\def\E{\mathcal E}

\newtheorem{thm}{Theorem}[section]
\newtheorem{prop}[thm]{Proposition}
\newtheorem{cor}[thm]{Corollary}
\newtheorem{lem}{Lemma}
\newtheorem*{prop*}{Proposition}
\newtheorem*{thm*}{Theorem}

\newtheorem*{claim*}{Claim}

\theoremstyle{definition}
\newtheorem*{df}{Definition}
\newtheorem*{dfs}{Definitions}

\theoremstyle{remark}
\newtheorem{rem}{Remark}
\newtheorem*{nota}{Notation}

\theoremstyle{remark}
\newtheorem{ex}{Example}

\DeclareMathOperator{\codim}{codim}

\DeclareMathOperator{\hh}{H}

\begin{document}

\title{On double Veronese embeddings in the Grassmannian $G(1,N)$}
\author{Jos\'e Carlos Sierra}
\address{Departamento de Algebra, Facultad de Matem\'aticas, Universidad Complutense, 
28040 Madrid}
\email{JoseCarlos\_Sierra@mat.ucm.es}
\author{Luca Ugaglia}
\address{Dipartimento di Matematica, Universit\`a degli Studi di Milano, Via Saldini 50,
20133 Milano}
\email{Luca.Ugaglia@unimi.it, }
\begin{abstract}
We classify all the embeddings of $\p^{n}$ in a Grassmannian of lines $G(1,N)$ 
such that the composition with Pl\"ucker is given by a linear system of 
quadrics of $\p^{n}$.
\end{abstract}
\maketitle

\section*{Introduction}

There exists a well known Hartshorne's Conjecture saying that a smooth $r$-di\-men\-sio\-nal
variety $X\su\p^{n}$ is a complete intersection if $n<\frac{3}{2}r$. A weaker statement
says that $X\su\p^{n}$ is a c.i. if $X$ has codimension 2 and $n\geq 7$ (see \cite{conjetura}). 
According to a theorem of Serre,
every codimension 2 smooth subvariety of $\p^{n}$, for $n\geq6$, can be given as the
zero locus of a section of a rank-2 vector bundle on $\p^{n}$. Therefore
the weaker Hartshorne's conjecture is equivalent to prove that every rank-2
bundle on $\p^{n}$ splits if $n\geq 7$ (indeed it is conjectured to be true for every
$n\geq 5$).
Moreover, to give a rank-2 bundle $\E$ on $\p^{n}$ together with an epimorphism 
$\oc_{\p^{n}}^{N+1}\to \E$ is equivalent to give a map
from $\p^{n}$ to a Grassmannian $G(1,N)$. If we tensor $\E$ with the line bundle
$\oc_{\p^{n}}(m)$, for a suitable $m>>0$, we get that the map associated
to the new rank-2 bundle is an embedding (see \cite{gh}). 
Therefore in order to prove Hartshorne's conjecture it is enough to consider bundles 
giving an embedding of $\p^{n}$ in $G(1,N)$.

When $\det(\E)=1$, there are only two embeddings of $\p^{n}$ in the
Grassmannian $G(1,N)$ such that the composition with the Pl\"ucker 
embedding of $G(1,N)$ gives rise to a linear space.
If $n\geq 3$, the only way is to take the rank-2 bundle 
$\oc_{\p^{n}}(1)\oplus\oc_{\p^{n}}$.

Let us consider a further step, i.e. rank-2 vector bundles $\E$ 
that give an embedding of $\p^{n}$ in a Grassmannian $G(1,N)$ and such that $\det(\E)=2$
(i.e. the composition with Pl\"ucker corresponds to the line bundle $\oc_{\p^{n}}(2)$).
In this paper we classify all such vector bundles obtaining as a corollary that
if $n\geq 4$ then $\E$ splits.

The paper is organized as follows. In the first section we recall some definitions and
results on embeddings in Grassmannians and on vector bundles on $\p^{n}$. 
In Section 2 we give some examples of double Veronese embeddings of $\p^{n}$ in 
$G(1,N)$ and we prove two lemmas (concerning the cases $n=1$ and $2$) that we will use
in Section 3 in order to prove our classification theorem.

We would like to thank Professor Fyodor Zak for his suggestion and his helpful discussions
on these topics and the organization of PRAGMATIC 2002 where this problem has been posed.

\section{Preliminaries}
\begin{nota}
Let us take an element of the Grassmannian $G(1,N)$; throughout the paper
we will use the following notations:\\
i) a small letter $l$, if we refer to it as a point of $G(1,N)$;\\
ii) a capital letter $L$, if we consider it as a subspace of $\p^{N}$.
\end{nota}

\begin{dfs}
A non-degenerate variety $X\su\p^{N}$ (i.e. not contained in a $\p^{N-1}$) is said to be 
{\em projective linearly normal} if $X$ is not 
projected from any non-degenerate variety contained in a bigger projective space, i.e. 
if $h^0(X,\oc_{X}(1))= N+1$.\\
A non-degenerate variety $X\su G(1,N)$ (i.e. not contained in a $G(1,N-1))$ is said to 
be {\em Grassmannian linearly normal} 
if $X$ is not projected from any non-degenerate variety contained in a bigger 
Grassmannian of lines, i.e. if
$h^0(X,\q_{|X})=N+1$ (where $\q$ is the universal quotient bundle of $G(1,N)$ and 
$\q_{|X}$ denotes its restriction $\q\otimes\oc_{X}$ to $X$). 
\end{dfs}
Let us recall some general facts about embeddings in Grassmannians of lines and Pl\"ucker 
embedding (for a detailed description see for instance \cite{trento}). \\ 
Giving a non-degenerate map $\varphi:X\rightarrow G(1,N)$
is equivalent to giving a rank-2 vector bundle $\E$ and an epimorphism 
$\phi:V\otimes\oc_{X}\rightarrow \E$, where
$V$ is an $(N+1)$-dimensional subspace of $\hh^0(X,\E)$. In this situation, 
$\E \cong \q_{|X}$. Moreover $\varphi$ is an embedding if different points of $X$ 
(maybe infinitely close) are mapped to different lines, i.e. if any subscheme of 
length two of $X$ imposes at least three conditions to $V$. 

Let us consider the embedding 
$X\stackrel{\bar{\varphi}}\hookrightarrow\p^{M}$, where $M={N+1\choose 2}-1$,
composition of the Pl\"ucker embedding and $\varphi$.
If the vector space $V$ is the whole $\hh^0(X,\E)$, $X$ is Grassmannian
linearly normal and $\bar{\varphi}$
is given by $M+1$ sections of the line bundle $\wedge^2\E$. 
Note that $X$ can be very degenerate in $\p^{M}$, since 
these sections are not necessarily independent, so we will always consider 
$X$ contained in its linear span $\langle X \rangle \cong \p^{r}\su\p^{M}$, 
where $r+1$ is the maximal number of 
independent sections of $\wedge^2\E$. In general, $X$ is not necessarily projective 
linearly normal in $\langle X \rangle$.
\begin{df}
Let us take a variety $X\su G(1,N)$, image of $\p^{n}$ via an embedding $\varphi$,
such that the composition $\bar{\varphi}$ with the Pl\"ucker embedding is a 
(maybe degenerate and not necessarily projectively linearly normal) double 
Veronese embedding $v_2(\p^{n})$ (i.e. $\wedge^2\E$ coincides with $\oc_{\p^{n}}(2)$). 
Throughout the paper we will say that $X$ is a {\em double Veronese embedding of 
$\p^{n}$ in $G(1,N)$}.
\end{df}

\vskip .2truecm

We are now going to recall some definitions and state some known results 
about vector bundles on $\p^{n}$ (for a detailed description see \cite{oko}).\\
Let $\E$ be a rank $r$ bundle on $\p^{n}$. According to a theorem of Grothendieck, 
for every $l\in G(1,n)$ there is an $r$-tuple
\[
a_{\E}(l)=(a_1(l),\ldots,a_r(l))\in\zz^r,
\]
with $a_1(l)\geq a_2(l)\geq\ldots\geq a_r(l)$, such that
\[
\E_{|L}=\E\otimes\oc_{L}=\bigoplus_{i=1}^r\oc_{\p^1}(a_i(l)).
\]
In this way can be defined a map
\[
a_{\E}:G(1,n)\rightarrow\zz^r.
\]
\begin{dfs}
The $r$-tuple $a_{\E}(l)$ is called the {\em splitting type of $\E$ on $L$}.\\
The bundle $\E$ is defined to be {\em uniform} if $a_{\E}$ is constant.\\ 
Let us now give $\zz^r$ the lexicographical order, i.e. $(a_1,\ldots,a_r)\leq
(b_1,\ldots,b_r)$ if the first non-zero difference $b_i-a_i$ is positive. We put
\[
{\bf a}_{\E}=\inf_{l\in G(1,n)}a_{\E}(l).
\]
The $r$-tuple ${\bf a}_{\E}$ is called the {\em generic splitting type of $\E$}.\\
A line $l\in G(1,n)$ is called a {\em jumping line} if $a_{\E}(l)>{\bf a}_{\E}$. The set of
jumping lines turns out to be a proper closed subset of the Grassmannian $G(1,n)$ (see
\cite{oko}).
\end{dfs}
\begin{thm}\label{rk=2}
Let $\E$ be a uniform rank 2 vector bundle on $\p^{n}$. Then either $\E$  
splits, or $n=2$ and $\E$ is a twist of the tangent bundle by some line bundle.
\end{thm}
\begin{proof}
See \cite{vdv}.
\end{proof}
\begin{thm}\label{line}
Let $\E$ be a rank $r$ vector bundle over $\p^{n}$, $x\in\p^{n}$ a point such that
$\E_{|L}=\oc_{\p^1}(a)^{\oplus r}$ for each line $L$ through $x$. Then $\E=\oc_{\p^{n}}(a)^{\oplus r}$. 
\end{thm}
\begin{proof}
We just have to apply \cite[Theorem 3.2.1]{oko} to the bundle $\E'=\E\otimes\oc_{\p^{n}}(-a)$. 
\end{proof}
\begin{thm}\label{plane}
A vector bundle $\E$ over $\p^{n}$ splits exactly when its restriction to
some plane $\Pi\su\p^{n}$ splits.
\end{thm}
\begin{proof}
See \cite[Theorem 2.3.2]{oko}.
\end{proof}

\section{Examples}
Let us now give some examples of Grassmannian linearly normal double Veronese embeddings 
$\p^{n}\stackrel{\varphi}\hookrightarrow G(1,N)$.
These will be the examples appearing in the statements of our main results.
\begin{ex}\label{o(1)+o(1)}
The rank-2 bundle $\E=\oc_{\p^{n}}(1)\oplus\oc_{\p^{n}}(1)$ gives an
embedding $X$ of $\p^{n}$ in the Grassmannian $G(1,2n+1)$ as the family of lines joining
the corresponding points on two disjoint $\p^{n}$'s. 
This is a double Veronese embedding, since $\wedge^2\E=\oc_{\p^{n}}(2)$.
\end{ex}
\begin{ex}\label{cone}
The rank 2 bundle $\E=\oc_{\p^{n}}(2)\oplus\oc_{\p^{n}}$ gives an
embedding $X$ of $\p^{n}$ in the Grassmannian $G(1,N)$, with 
$N={n+2\choose 2}$. This is the family of ruling lines of a cone 
over the double Veronese embedding $v_2(\p^{n})\su\p^{N-1}$, with vertex
a point. Again $\wedge^2\E=\oc_{\p^{n}}(2)$, and hence $X$ is a double 
Veronese embedding of $\p^{n}$.
\end{ex}
\begin{ex}\label{bisec}
The family of the bisecant lines to a rational normal cubic 
is a double Veronese embedding of $\p^2$ in $G(1,3)$ (see \cite{arr1}). 
In this case the bundle $\E$ is a Steiner bundle, i.e. it is the dual
of the kernel of a map $\oc_{\p^2}^{\oplus 4}\rightarrow\oc_{\p^2}(1)\oplus\oc_{\p^2}(1)$, 
corresponding to the choice of $4$ general sections of the bundle $\oc_{\p^2}(1)\oplus\oc_{\p^2}(1)$.
From the exact sequence 
\[
\xymatrix{ 0 \ar[r] & \oc_{\p^2}(-1)\oplus\oc_{\p^2}(-1)\ar[r] & \oc_{\p^2}^{\oplus 4}\ar[r] & 
\E\ar[r] & 0,}
\]
we get that $h^0(\p^2,\E)=4$, and hence the surface is Grassmannian linearly normal.
\end{ex}
\begin{ex}\label{quadric}
The family $X$ of lines contained in a smooth hyperquadric $Q\su\p^4$ is a double Veronese 
embedding of $\p^3$ in $G(1,4)$ (see \cite{tan}). The vector bundle $\E$ is the cokernel of
a map $\oc_{\p^3}\rightarrow\Omega_{\p^3}(2)$ corresponding to the choice of a general section
of the twist of the cotangent bundle $\Omega_{\p^3}$ by $\oc_{\p^3}(2)$. 
From the sequence
\[
\xymatrix{ 0 \ar[r] & \oc_{\p^3}\ar[r] & \Omega_{\p^3}(2)\ar[r] & \E\ar[r] & 0,}
\]
and the Euler sequence of $\p^3$ we get that $h^0(\p^3,\E)=5$ and
hence $X$ is Grassmannian linearly normal in $G(1,4)$. 
\end{ex}
\begin{ex}\label{2,2}
Taking the restriction of the embedding above to a general plane in $\p^3$ we
get a double Veronese embedding of $\p^2$ in $G(1,4)$.
Geometrically speaking it is the set of lines contained in $Q$ and meeting a 
fixed line in it. It is again Grassmannian linearly normal by the same reason of 
Example \ref{quadric}. 
\end{ex}
\begin{rem}\label{proje}
We recall that the variety given by the bundle $\oc_{\p^{n}}(1)\oplus\oc_{\p^{n}}(1)$ can be 
isomorphically projected from $G(1,2n+1)$ to $G(1,m)$, for $n+1\leq m\leq 2n$ 
(see \cite{arr2}). The variety given by $\oc_{\p^{n}}(2)\oplus\oc_{\p^{n}}$ can be 
isomorphically projected from $G(1,N)$ to $G(1,m)$, for $2n+1\leq m\leq N-1$, 
since the Veronese variety $v_2(\p^{n})$ can be projected from $\p^{N-1}$ to $\p^{m-1}$.\\
Conversely, varieties of Examples \ref{bisec}, \ref{quadric} and \ref{2,2} 
cannot be isomorphically projected to a smaller Grassmannian. This claim is 
obvious for Example \ref{bisec}, while for Examples \ref{quadric} and \ref{2,2} 
is enough to realise that through the general point of $\p^4$ there pass a 
2-dimensional family of planes intersecting $Q$ along a degenerate conic, 
i.e. two lines. These two lines are projected to the same line, giving rise 
to a singularity.
\end{rem}
We end the section with the classification of double Veronese 
embeddings of $\p^1$ and $\p^2$ respectively. This will be useful to simplify the 
proof of Theorem \ref{class}. 
\begin{lem}\label{p1}
The only Grassmannian linearly normal double Veronese embeddings of $\p^1$ in a 
Grassmannian of lines are as in Examples \ref{o(1)+o(1)} or \ref{cone}.
\end{lem}
\begin{proof}
Let us denote by $X$ the image of $\p^1$ via the double Veronese embedding
$\varphi:\p^1\hookrightarrow G(1,N)$,
corresponding to the bundle $\E=\oc_{\p^1}(\alpha)\oplus\oc_{\p^1}(\beta)$. 
Since $\varphi$ is an embedding,
we must have $\alpha,\beta\geq 0$. Moreover, by definition  it must be
$\wedge^2\E=\oc_{\p^1}(2)$, which implies $\alpha+\beta=2$. Therefore the only two
possibilities are either $\E=\oc_{\p^1}(1)\oplus\oc_{\p^1}(1)$ (corresponding to 
Example \ref{o(1)+o(1)}) or $\E=\oc_{\p^1}(2)\oplus\oc_{\p^1}$ (corresponding to 
Example \ref{cone}). We remark that in the first case $X$ is one of the 
rulings of a smooth quadric surface (and can be projected to $G(1,2)$ as 
the tangent lines to a smooth conic). In the second case $X$ corresponds 
to the ruling lines of a cone over a smooth conic.
\end{proof}
\begin{dfs}
Given a surface $S\su G(1,N)$, in the Chow ring of $G(1,N)$ we can write 
$[S]=a\Omega(0,3)+b\Omega(1,2)$, where $a=[S]\cdot\Omega(N-3,N)$ is the {\em order} of $S$
and $b=[S]\cdot\Omega(N-2,N-1)$ is its {\em class}.\\
The pair $(a,b)$ is defined to be the {\em bidegree} of $S$.\\
The {\em degree} of $S$ via Pl\"ucker embedding turns out to be $a+b$. 
\end{dfs}
\begin{prop}\label{p2}
The only Grassmannian linearly normal double Veronese embeddings of $\p^2$ in a 
Grassmannian of lines are as in Examples \ref{o(1)+o(1)}, \ref{cone}, 
\ref{bisec} or \ref{2,2}.
\end{prop}
\begin{proof} 
Let us denote by $\E$ the rank-2 bundle giving the embedding of $X=\p^2$ in a 
Grassmannian of lines $G(1,N)$, i.e. $\E \cong \q_{|X}$, and consider its 
restriction $\E_{|L}=\E\otimes\oc_{L}$ to a general line $L\su\p^2$. 
We remark that $\E_{|L}$ gives a double Veronese embedding of $L$ since it is the 
restriction of the embedding given by $\E$. Hence by Lemma \ref{p1}
either $\E_{|L}=\oc_{\p^1}(2)\oplus\oc_{\p^1}$ or $\E_{|L}=\oc_{\p^1}(1)\oplus\oc_{\p^1}(1)$, 
corresponding to ${\bf a}_{\E}=(2,0)$ or $(1,1)$ respectively.

\vskip .1truecm

If the generic splitting type ${\bf a}_{\E}$ is $(2,0)$, then there are no jumping lines 
and $\E$ must be uniform. Hence, by Theorem \ref{rk=2}, either $\E$ splits, or $\E$ is the twist by a line bundle of the tangent bundle. 
But this last possibility cannot occur, since in this case the line bundle 
$\wedge^2\E$ would have odd degree and hence
it could not give the double Veronese embedding of $\p^2$.
Therefore, since $\E$ is decomposable and ${\bf a}_{\E}=(2,0)$, it must be
$\E=\oc_{\p^2}(2)\oplus\oc_{\p^2}$. Moreover we assume $X$ Grassmannian linerly normal, 
so we take $V$ to be the whole $\hh^0(\p^2,\E)$, 
corresponding to the double Veronese embedding of Example \ref{cone} for $n=2$.

\vskip .1truecm

If ${\bf a}_{\E}=(1,1)$, for a jumping line $L$ we must have 
$\E_{|L}=\oc_{\p^1}(2)\oplus\oc_{\p^1}$, which is equivalent to say that $L$ is embedded 
in a $G(1,3)\subseteq G(1,N)$ as the rulings of a quadratic cone.
Let us denote by $\j\su G(1,2)$ the closed set of jumping lines of $\E$. 
We distinguish two different cases depending on the codimension of $\j$.
\begin{enumerate}
\item[{\bf i)}] $\codim\j\geq 2$. 
In this case 
there are at most finitely many jumping lines for $\E$. In particular, 
through a general point $x\in\p^2$ there pass no jumping lines, which means that
$\E_{|L}=\oc_{\p^1}(1)\oplus\oc_{\p^1}(1)$ for each $L$ through $x$. By Theorem \ref{line}
we get that $\E=\oc_{\p^2}(1)\oplus\oc_{\p^2}(1)$ and hence $\E$ 
gives the double Veronese embedding of Example \ref{o(1)+o(1)} 
when we consider all the sections of $\hh^0(\p^2,\E)$.
\item[{\bf ii)}] $\codim\j=1$. 
For every irreducible maximal component of $\j\su G(1,2)$, there exists a 
fundamental curve $C\su \p^{N}$, i.e. a curve which is cut by all the lines of the 
surface. The points of $C$ are the vertices of the quadratic cones.\\
The number of fundamental points of $C$ contained in a general line of $X$ cannot 
be bigger than two by the classical {\em trisecant lemma}, so there are just two 
possibilities: either $\p^2$ is embedded as the family of bisecant lines to the 
curve $C\su\p^{N}$, or $C$ is a line of $X$ contained in all the quadratic cones 
corresponding to the component of $\j$ cited above. \\
In the first case, the bisecant lines 
passing through a general point $c\in C$ give rise to a quadratic cone 
(since they correspond to the embedding of a jumping line). This implies that $C$ 
is contained in a $\p^3$ and that its projection from a general point on it is a 
smooth conic, and hence $C$ must be a rational normal cubic of $\p^3$. In this 
way $X$ turns out to be the 
Veronese surface of Example \ref{bisec}. Note that $\j\subset G(1,2)$ is a conic 
embedded as the family of tangent lines to a smooth conic of $\p^2$. This conic of 
$\p^2$ is embedded in $G(1,3)$ as the tangent developable to the rational normal cubic $C$.\\
If $C$ is a line we claim that the bidegree $(a,b)$ of the surface $X$ 
must be $(2,2)$. In fact, the class $b$ of 
$X$ is the number of lines contained in a general 
hyperplane $H$ of $\p^{N}$. But $H$ contains exactly two lines, 
corresponding to the hyperplane section of the quadratic cone passing through $H\cap R$
and hence $b=2$. Moreover, since $X$ is a Veronese surface, its degree is 
$a+b=4$, which proves our claim. This also implies that the $3$-fold $\overline X$ 
covered by the lines of $X$ is either a hyperquadric $Q\subset\p^4$ or a $\p^3$.  The later case is 
not possible since there are no Veronese surfaces of bidegree $(2,2)$ in $G(1,3)$ 
(see \cite{arr1}). \\
Let us see that $Q$ is smooth. If $Q$ is a quadric of rank $2$ or $3$, then it contains 
a family of planes. Since $a=2$ through the general point of $Q$ there pass just 
one line of $X$, so the lines of $X$ on such planes move on a pencil through a point. 
These pencils are embedded as lines via Pl\"ucker, which is absurd because the Veronese 
surface does not contain lines.\\
Therefore $X$ is the set of lines of a smooth quadric $Q$ meeting a line contained in it,  
which is the Veronese surface of Example \ref{2,2}. 
We remark that in this case $\j\subset G(1,2)$ is a pencil of lines. 
\end{enumerate}
\end{proof}

\section{Classification}
In this section we classify all double Veronese embeddings $X$ of $\p^{n}$ in $G(1,N)$.
In order to do that, we first consider Grassmannian linearly normal
double Veronese embeddings of $\p^{n}$.
\begin{thm}\label{class}
Varieties of Examples \ref{o(1)+o(1)}, \ref{cone}, \ref{bisec}, \ref{quadric} and 
\ref{2,2} are the only smooth, Grass\-man\-nian linearly normal, double Veronese embeddings of 
$\p^{n}$ in $G(1,N)$.
\end{thm}
\begin{proof}
Let us denote by $\E$ the rank-2 bundle giving the embedding of $\p^{n}$ in $G(1,N)$.
Arguing as we did in the proof of Proposition \ref{p2} we can say that
${\bf a}_{\E}=(2,0)$ or $(1,1)$. Moreover, if ${\bf a}_{\E}$ is $(2,0)$, 
we can conclude as before that $X$ is the double Veronese embedding of Example \ref{cone}.
If ${\bf a}_{\E}=(1,1)$, for a jumping line $L$ we must have $\E_{|L}=\oc_{\p^1}(2)\oplus\oc_{\p^1}$.
We still denote by $\j\su G(1,n)$ the closed set of jumping lines of $\E$ and
we distinguish two cases.
\begin{enumerate}
\item[{\bf i)}] $\codim\j\geq 2$. 
Let us take a plane $\Pi\su\p^{n}$ and consider the restriction 
$\E'=\E\otimes\oc_{\Pi}$. This gives
an embedding of $\Pi$ in a Grassmannian of lines as a Veronese surface.  Since there are 
at most a finite number of jumping lines on $\Pi$, it follows from Proposition \ref{p2} 
that $\E'$ splits.
By Theorem \ref{plane} also $\E$ splits. We conclude that $\E=\oc_{\p^{n}}(1)\oplus\oc_{\p^{n}}(1)$, 
which gives the double Veronese embedding of Example \ref{o(1)+o(1)} when we consider 
all the sections of $\hh^0(\p^{n},\E)$.
\item[{\bf ii)}] $\codim\j=1$. 
Let us consider the following incidence diagram:
\[
\xymatrix{
& I\ar[dr]_{p_2} \ar[dl]^{p_1}\\
X && \j,}
\]
where we put $I=\{(l,m)\in X\times\j\mid  l\in M\}$. The general fiber of $p_2$ has
dimension 1, and hence $\dim I=2n-2$ and $\dim p_1^{-1}(l)=n-2$, for a general $l\in X$.
This is equivalent to say that the general line $L$ 
is contained in an $(n-2)$-dimensional family of cones whose ruling lines correspond 
also to points of $X$.
In particular the general line $L$ meets an $(n-1)$-dimensional family of lines of $X$.
Therefore, either there exists one point through which there pass
an $(n-1)$-dimensional family of lines, or through the general point of $L$
there pass an $(n-2)$-dimensional family of lines of $X$ (with $n\geq 3$).\\ 
In the first case we get that there exists a fundamental curve $C\su\p^{N}$ such that through
the general point $c\in C$ there pass an $(n-1)$-dimensional family of lines of $X$.
We expect that through two general points $c_1,c_2\in C$ there
pass an $(n-2)$-dimensional family of lines, which implies $n-2=0$. 
Therefore $X\cong\p^2$ and the classification follows from Proposition \ref{p2}.\\
Let us consider now the second case, i.e. through a general point of $L$ there pass an
$(n-2)$-dimensional family of lines of $X$. We denote by $\overline X\su\p^{N}$ the union of 
the lines of $X$ and consider the incidence diagram
\[
\xymatrix{
& W\ar[dr]_{q_2} \ar[dl]^{q_1}\\
X && \overline X,}
\]
where we put $W=\{(l,y)\mid y\in L\}$. Looking at the first projection we get $\dim(W)=n+1$
and, since through a general point $y\in \overline X$ there pass an $(n-2)$-dimensional 
family of lines
of $X$, $\dim \overline X=n+1-(n-2)=3$. In particular $\overline X$ is a 3-dimensional 
projective variety containing
a 3-dimensional family of lines (if $\overline X$ contains a bigger family of lines 
then $\overline X=\p^3$ and $X\su G(1,3)$, but this is not possible since $v_2(\p^3)$ 
cannot be embedded in $\p^5$). Then either $\overline X$ is swept out by a 1-dimensional
family of planes, or it is a hyperquadric of $\p^4$. The former is not possible since in this case
the 3-dimensional family of lines contained in $\overline X$ would be a scroll of planes. 
Conversely, the later corresponds to the double Veronese embedding of Example \ref{quadric}. 
Note that $\j\subset G(1,3)$ is a hyperplane section of $G(1,3)$. Moreover, jumping lines 
$L\in\j$ correspond to tangent spaces to $Q$, since they intersect $Q$ along quadric cones. 
The rulings of the cone give the embedding of $L$.
\end{enumerate}
\end{proof}
\begin{rem}
All the double Veronese embeddings we classified above are 
projective linearly normal in $\langle X \rangle$ via Pl\"ucker embedding. 
This is obvious for Examples \ref{o(1)+o(1)}, \ref{cone} and \ref{bisec}. For 
Example \ref{quadric} see \cite{tan} and, since Example \ref{2,2} is the restriction to a 
plane of Example \ref{quadric} it is also projective linearly normal. 
\end{rem}
In order to complete the classification we are now going to consider
double Veronese embeddings that can be projected to a smaller Grassmannian.
In particular we study when the Grassmannian projection gives rise to a 
variety which is not projective linearly normal in $\langle X \rangle$.
In this way we also classify embeddings of $\p^{n}$ in a $G(1,N)$ such that
the composition $\bar{\varphi}$ with Pl\"ucker is given by a proper subspace  of 
$\hh^0(\p^{n},\oc_{\p^{n}}(2))$.
\begin{prop}
If $X\su G(1,N)$ is a double Veronese embedding of
$\p^{n}$ and it is not projective linearly normal, then $X$ is a projection 
of the variety of Example \ref{cone}.
\end{prop}
\begin{proof}
By Remark \ref{proje}, varieties of Examples \ref{o(1)+o(1)} and \ref{cone} are the 
only double Veronese embeddings that 
can be isomorphically projected to a smaller Grassmannian of lines.
In order to prove the proposition we show that when we project, only in the first
case we can obtain a non-projective linearly normal variety.

Let us take the embedding of $\p^{n}$ in the Grassmannian 
$G(1,N)$, with $N={n+2\choose 2}$, given by the bundle $\oc_{\p^{n}}(2)\oplus\oc_{\p^{n}}$.
The image $X\su G(1,N)$ can be described as the family of ruling lines of a cone 
over the double Veronese embedding 
$v_2(\p^{n})\su H\cong\p^{N-1}$, with vertex a point $v\notin H$. 
Let us take a linear space $L\su H$, of dimension $k$, with 
$0\leq k\leq \frac{(n+1)(n-2)}{2}$, which does not intersect the secant variety of $v_2(\p^{n})$.
We put $m=N-k-1$. The projection $\pi_{L}:G(1,N)\rightarrow G(1,m)$, restricted to $X$, 
is an isomorphism.
The image $X'$ can be described as the family of ruling lines of a cone 
over the projection of $v_2(\p^{n})$ to $\p^{m-1}$. Therefore the composition $\bar{\varphi}$
with Pl\"ucker embedding is given by a  subspace of $\hh^0(\p^{n},\oc_{\p^{n}}(2))$ of dimension 
$m$ and $X'$ is not projective linearly normal.

Let us take now the embedding $X$ of $\p^{n}$ in $G(1,2n+1)$
given by the bundle $\oc_{\p^{n}}(1)\oplus\oc_{\p^{n}}(1)$. We prove that if we project $X$ to 
$G(1,m)$, with $n+1\leq m\leq 2n$ and we consider the composition with Pl\"ucker embedding, 
we always get the map associated to the complete linear system of quadrics $|\oc_{\p^{n}}(2)|$.
It is enough to prove it for the projection to $G(1,n+1)$. We recall that $X$ can be seen
in $G(1,2n+1)$ as the family of lines joining corresponding points on two disjoint 
$\p^{n}$'s. The same geometric description holds after projecting to $G(1,n+1)$, but here the two
$\p^{n}$'s, say $V_0$ and $W_0$, intersect in a $\p^{n-1}$. Let us fix coordinates
$(x_0:\ldots:x_{n+1})$ for $\p^{n+1}$ and let us denote by 
$\phi:V_0\rightarrow W_0$ the correspondence.
By induction we construct the linear spaces
$W_i:=\phi(V_i)$ and $V_{i+1}:=V_i\cap W_i$, for $i=0,\ldots,n$. Since the projection 
is an isomorphism, we have that no line of $X$ is contracted, and this implies
that $V_{i+1}\subsetneq V_{i}$, or $\dim(V_i)=\dim(W_i)=n-i$.
Changing coordinates we can suppose that $V_0=\{x_{n+1}=0\}$,  
$V_i=\{x_{n+1}=x_0=\ldots=x_{i-1}=0\}$ for $i=1,\ldots,n$ and 
$W_i=\{x_0=\ldots=x_{i}=0\}$. Under these assumptions $X$ can be described as
the family of lines spanned by the rows of the matrix
\begin{equation}\label{coord}
\left (
\begin{array}{ccccc}
t_0 & t_1 & \cdots & t_n & 0 \\
0 & l_0 & \cdots & l_{n-1} & l_n
\end{array}
\right ),
\end{equation}
where $t_0,t_1,\ldots,t_n$ are homogeneous coordinates of $\p^{n}$ and 
$l_i=a_{i,i}t_i+a_{i,i+1}t_{i+1}+\ldots a_{i,n}t_n$ is a linear 
form involving only the last $n-i-1$ variables.
We can view $\phi$ as the map from $\p^{n}$ to $\p^{n}$, sending $t_i$ to $l_i$, 
and hence it is represented by a lower triangular matrix $T$,
whose determinant is $\prod_{i=0}^na_{i,i}$. We remark that, since $\phi$ is an isomorphism, 
$a_{i,i}\neq 0$ for $i=0,\ldots,n$. 
Let us substitute to $x_1$ a suitable linear combination of $x_1,\ldots,x_{n+1}$, in order to 
get $t_0$ in the second place of the second row of the matrix above. The corresponding 
element on the first row is now a linear combination $t_1'$ of $t_1,\ldots,t_n$.
Let us write now $l_1$ with respect to $t_1',t_2,\ldots,t_n$. As before we can send 
$x_2$ to a suitable combination of $x_2,\ldots,x_{n+1}$ in order to get $t'_1$ 
in the third place of the second row of the matrix \ref{coord}, and so on. 
In this way we get a base change of $\p^{n+1}$ and of $V_0$, since
the corresponding matrices are triangular and the elements on the diagonal
are products of some $a_{i,i}^{-1}$. With the new bases the variety $X$ can be described by
\[
\left (
\begin{array}{ccccc}
t_0 & t'_1 & \cdots & t'_n & 0\\
0 & t_0 & \cdots & t'_{n-1} & t'_n
\end{array}
\right ),
\]
whose minors give a basis for the space of degree $2$ polynomials in $t_0,t'_1,\ldots,t'_n$.
\end{proof}

\begin{rem}
In \cite{tan}, H. Tango classified embeddings of $\p^{n}$ in $G(1,n+1)$. 
There are just 4 possibilities, namely, the {\em star of lines}, Examples 
\ref{bisec} and \ref{quadric}, and the projection to $G(1,n+1)$ of Example 
\ref{o(1)+o(1)}. As a corollary of our classification we get that {\em all} 
double Veronese embeddings except the cone case of Example \ref{cone} and the Veronese surface of Example \ref{2,2} fit in a $G(1,n+1)$.
\end{rem}

Finally, we state the result in connection with Hartshorne's 
Conjecture quoted in the introduction. 
\begin{cor}\label{split}
Rank-$2$ vector bundles over $\p^{n}$ giving a double Veronese embedding split if $n\geq 4$.
\end{cor}

\providecommand{\bysame}{\leavevmode\hbox to3em{\hrulefill}\thinspace}

\end{document}